\documentclass[11pt]{amsart}

\usepackage{amscd, amssymb}
\usepackage{enumerate}
\usepackage[all]{xypic}
\usepackage{multirow}
\usepackage{hhline}
\usepackage{verbatim}
\usepackage{mathdots}
\usepackage{comment}
\usepackage{tikz-cd}
\usepackage[titletoc]{appendix}
\usepackage{stmaryrd}%this conflicts with extpfeil
\usepackage{graphicx}
\usepackage{mathtools}
\usepackage{leftidx}
\usepackage{mathabx }
\usepackage{amsfonts}
\usepackage{mathrsfs}

\usepackage{dynkin-diagrams}
\usepackage{xcolor}

\usepackage{enumitem}

\usepackage[margin=0.8in]{geometry}

\usepackage[colorlinks=true, linkcolor=blue,
 citecolor=blue, urlcolor=blue]{hyperref}

\theoremstyle{plain}
\newtheorem{Thm}[equation]{Theorem}
\newtheorem*{Thm*}{Theorem}
\newtheorem{Cor}[equation]{Corollary}
\newtheorem{Prop}[equation]{Proposition}
\newtheorem{Lem}[equation]{Lemma}

\theoremstyle{definition}
\newtheorem*{Rmk}{Remark}

\numberwithin{equation}{section}

\newcommand{\C}{\mathbb{C}}

\newcommand{\Q}{\mathbb{Q}}

\newcommand{\Z}{\mathbb{Z}}

\newcommand{\Gt}{\overline{G}}

\newcommand{\It}{\overline{I}}

\newcommand{\Kt}{\overline{K}}
\newcommand{\Tt}{\overline{T}}

\newcommand{\Yt}{Y_{Q,n}}

\newcommand{\Gbb}{{\mathbb{G}}}

\newcommand{\ttt}{\overline{t}}

\newcommand{\yt}{\overline{y}}

\newcommand{\s}{\mathbf{s}}

\newcommand{\Ocal}{\mathcal{O}}

\newcommand{\Hcal}{\mathcal{H}}

\newcommand{\Kf}{\mathsf{K}}

\newcommand{\la}{\langle}
\newcommand{\ra}{\rangle}

\title[High degree covers and genericity]{High degree covers and genericity}

\author{Edmund Karasiewicz}


\allowdisplaybreaks

\begin{document}

\begin{abstract}
    It is known that the genuine irreducible Iwahori-spherical representations of a tame $n$-fold Brylinski-Deligne cover of a $p$-adic split connected reductive group are all generic, when the degree $n$ is suitably large. This has been established by studying the the structure of the Gelfand-Graev representation as a module under the action of the Iwahori-Hecke algebra.  
    
    In this paper, we present a different proof based on metaplectic Moy-Prasad theory. 
\end{abstract}

\maketitle

\setcounter{tocdepth}{1}
%\tableofcontents

%%%%%%%%%%%%%%%%%%%%%%%%%%%%%%%%%%%%%%%%%%%%%%%%%%%%%%%%%%%%%%%%%%%%%%%

\section{Introduction}

%%%%%%%%%%%%%%%%%%%%%%%%%%%%%%%%%%%%%%%%%%%%%%%%%%%%%%%%%%%%%%%%%%%%%%%  

In the representation theory of tame $n$-fold Brylinski--Deligne (BD) covers of split connected reductive $p$-adic groups, it is known that genuine irreducible representation in the Iwahori-spherical block (equivalently the block of unramified principal series \cite{S04}) are all generic, if $n$ is suitably large. This was proved in \cite{GGK24,GGK25} by examining the Gelfand-Graev representation as a module over a Hecke algebra. Earlier work on the genericity of theta representations (also called exceptional representations) appears in Kazhdan--Patterson \cite{KP84} for Kazhdan-Patterson covers and Gao \cite{Gao_theta_17} for more general BD covers.

The objective of this note is to present a simple explanation for this phenomenon based on Moy--Prasad theory \cite{MP94,MP96}, which was extended to tame covers by Howard--Weissman \cite{HW09}. 

We note that Moy--Prasad theory already plays a role in recent calculations of the stable wave front set of certain representations in the Iwahori-block of tame BD covers \cite{KOW26,Gao-Wang26}. However, these works involve the local character expansion and as a result impose constrains on the residual characteristic beyond the tame assumption, $p$ does not divide $n$. Here, we only assume that $p$ does not divide $n$.

%The more general question of computing the (stable) wave front set of representations in the Iwahori-block was recently considered in \cite{Gao-Wang26}.

%(See \cite{Gao_theta_17} for an extension to BD covers.) The result in the more general setting above is proved in \cite{GGK24} and further examined in the follow-up \cite{GGK25}. 

%See also \cite{Gao-Wang26}, for calculations of the wave front set of certain Iwahori-spherical representations.

Let $F$ be a finite extension of $\Q_{p}$. Let $\Ocal\subset F$ be the ring of integers and fix $\varpi\in\Ocal$ a uniformizer. Let $\kappa$ be the residue field of $F$ and $q$ its cardinality. 

For a field $L$ and $m\in\Z_{\geq 1}$, we will write $\mu_{m}(L)$ for the group of $m$-th roots of unity in $L$. For the $p$-adic field $F$ we set $\mu_{m}:=\mu_{m}(F)$. 

%Let $\Gbb$ be a split connected reductive linear algebraic group over $F$. 

Let $G$ be the group of $F$-points of a split connected reductive linear algebraic group over $F$. Let $T\subset G$ be a maximal split torus with character lattice $X$ and cocharacter lattice $Y$. Let $\Phi\subset X$ be the subset of roots and $\Phi^{\vee}\subset Y$ the subset of corroots attached to the pair $(G,T)$. For $\alpha\in\Phi$ we write $U_{\alpha}$ for the root subgroup attached to $\alpha$.

Let $B=B^{+}\supset T$ be a Borel subgroup with unipotent radical $U$. Let $\Delta\subset\Phi$ be the set of simple roots determined by $B$. Let $B^{-}\supset U^{-}$ be the opposite Borel and unipotent radical.

Let $K$ be a hyperspecial maximal compact subgroup such that $T_{0}:=T\cap K$ is the maximal compact subgroup of $T$. Let $I^{\pm}\subset K$ be an Iwahori subgroup such that $B^{\pm}\cap K = B^{\pm}\cap I$. For a pro-finite group $H$, let $H_{+}$ be its maximum pro-$p$ subgroup, when it exists. We will write $I^{-}_{+}:=(I^{-})_{+}$.

Assume that $|\mu_{n}|=n$. Let $\Gt$ be a topological central extension of $G$ by $\mu_{n}$ arising from the construction of Brylinski-Deligne \cite{BD01}. The group $\Gt$ fits into a central extension of the form
\begin{equation}\label{Eqn:CExt}
    \mu_{n}\hookrightarrow \Gt\stackrel{\mathrm{pr}}{\twoheadrightarrow} G.
\end{equation}
For any $H\subset G$, let $\overline{H}:=\mathrm{pr}^{-1}(H)$. Let $\epsilon=\epsilon_{n}:\mu_{n}\rightarrow \C^{\times}$ be a faithful character. A $\Gt$-representation is called $\epsilon$-genuine if $\mu_{n}$ acts through $\epsilon$. Let $\mathrm{Irr}_{\epsilon}(\Gt)$ be the set of isomorphism classes of $\epsilon$-genuine irreducible smooth $\Gt$-representations.

Recall that the covering group $\Gt$ comes with $Q:Y\rightarrow \Z$ a $W$-invariant quadratic form, where $W$ is the Weyl group of $T$. Let $B=B_{Q}$ be the bilinear form on $Y$ associated to the quadratic form $Q$. For any $m\in\Z$, let $Y_{Q,m}:=\{y\in Y|B(y,z)\in m\Z\text{, for all }z\in Y\}$.

For the remainder of this note we fix $n\in\Z_{\geq 1}$ such that $\gcd(p,n)=1$ and $|\mu_{n}|=n$. This implies that $n$ divides $q-1$. In particular, $\Gt$ is a tame BD cover.

We will also assume that there exists $\s_{K}:K\rightarrow \Kt$, a splitting of the sequence \eqref{Eqn:CExt} over $K$. This happens, for example, when $\Gt$ is incarnated by $(D,\eta)$, where $D$ is a fair bisector and the group homomorphism $\eta:\Z[\Phi^{\vee}]\rightarrow F^{\times}$ is valued in $\Ocal^{\times}F^{\times n}$ \cite[Theorem 4.2 and subsequent remarks]{GG18}.

Now we state our main theorem.

\begin{Thm}\label{Thm:Main_Thm}
    Let $\Gt$ be an $n$-fold BD cover of a split connected reductive group $G$ over a $p$-adic field $F$ such that $|\mu_{n}|=n$, $\gcd(p,n)=1$, and $K$ splits \eqref{Eqn:CExt} with splitting $\s_{K}$. Let $\pi\in \mathrm{Irr}_{\epsilon}(\Gt)$. Let $N$ be the number of free $W$-orbits in $Y/\Yt$. 

    If $\pi^{\s_{K}(I)}\neq 0$, then the Whittaker dimension with respect to any non-degenerate $\psi\in \mathrm{Hom}(U^{-},\C^{\times})$ of conductor $(\varpi)$ is at least $N$. 
\end{Thm}

\begin{Rmk}
    Earlier we used the informal statement `$n$ is suitably large' (this also hides an implicit dependence on $Q$). In the theorem this informal statement is replaced by `$Y/\Yt$ contains a free $W$-orbit'. 
\end{Rmk}

Actually, we will prove the more general Theorem \ref{Thm:generalized_main}, from which Theorem \ref{Thm:Main_Thm} follows as a corollary.

Our main theorem is the result of two inputs. The first input follows from unraveling some results from depth $0$ metaplectic type theory. The main consequence of this is that in the setting of Theorem \ref{Thm:Main_Thm}, the $\Gt$-representations contain regular principal series among their $\Kf$-types, which are generic. The second input is a result lifting Whittaker functionals from finite groups to $p$-adic groups and their covers.

\section{Metaplectic Moy-Prasad theory}

Moy-Prasad theory was developed in \cite{MP94,MP94}, in the setting of linear groups. These ideas were generalized to tame covering groups by Howard--Weissman \cite{HW09}. We will specialize some results of Moy--Prasad and Howard--Weissman and then examine their effect on the existence of Whittaker functionals. 

The main point for this paper is that the set of genuine depth $0$ Iwahori types admits a free action of the finite abelian group $Y/\Yt$, and this proliferation of distinct types can create many independent Whittaker functionals when the degree of the cover $n$ is suitably large. 

We begin with the commutator of $\Tt$, which is a source of non-genuine characters. For $y\in Y$, let $\phi_{y}\in \mathrm{Hom}(T,\mu_{n}(\C))$ be the character defined as follows. Let $t\in T$ and $y\in Y$. Fix a uniformizer $\varpi\in F^{\times}$. Let $\ttt\in \mathrm{pr}^{-1}(t)$ and $\yt(\varpi)\in \mathrm{pr}^{-1}(y(\varpi))$ be arbitrary. Let $[-,-]$ be the commutator in $\Tt$. We define
\begin{equation*}
    \phi_{y}(t):= \epsilon([\ttt,\yt(\varpi)]).
\end{equation*}
Under the identification $T\cong Y\otimes F^{\times}$, for $z\in Y$, $a\in F^{\times}$ we have
\begin{equation}\label{Eq:Tt_Comm_Form}
    \phi_{y}(z\otimes a) = (\varpi,a)_{n}^{B(y,z)},
\end{equation}
where $(-,-)_{n}$ is the $n$-th order Hilbert symbol for $F$.

The next lemma follows from well-known properties of the tame symbol $(-,-)_{n}$ and identity \eqref{Eq:Tt_Comm_Form}. 

\begin{Lem}\label{Lem:phi_y_props}
    Assume that $p$ does not divide $n$. Then:
    \begin{enumerate}
        \item\label{Lem:phi_y_props:1} $\phi_{y}|_{T_{0}}$ does not depend on the choice of $\varpi$;
        \item\label{Lem:phi_y_props:2} $\phi_{y}$ is trivial on $T_{1}:=T\cap I_{+}$, the maximal pro-$p$ subgroup of $T$.
        \item\label{Lem:phi_y_props:3} $\phi_{y}|_{T_{0}}$ is trivial if and only if $y\in \Yt$.
    \end{enumerate}
\end{Lem}

By Lemma \ref{Lem:phi_y_props} Item \eqref{Lem:phi_y_props:2}, $\phi_{y}|_{T_{0}}$ factors through $T_{\kappa}:=T_{0}/T_{1}$. 

\begin{Lem}\label{Lem:Char_equi}
    The map $Y/\Yt\hookrightarrow\mathrm{Hom}(T_{\kappa},\mu_{n}(\C))$ defined by $y\mapsto \phi_{y}|_{T_{0}}$ is a $W$-equivariant injective group homomorphism.

    %Moreover, if $\Yt=nY$, then the map is a bijection.
\end{Lem}

\begin{proof}
    One can directly check that the map is a $W$-equivariant group homomorphism. By Lemma \ref{Lem:phi_y_props} Item \eqref{Lem:phi_y_props:3}, the character $\phi_{y}$ is trivial if and only if $y\in \Yt$. Thus the map is injective.
\end{proof}

The inclusion $T_{0}\subset I$ induces an isomorphism $T_{\kappa}\cong I/I_{+}$. Thus we can also view $\phi_{y}|_{T_{0}}$ as a character of $I$. From now on we will just write $\phi_{y}$, instead of $\phi_{y}|_{T_{0}}$, for this character on $T_{0}$, $T_{\kappa}$, or $I$.

Let $\s_{I}:=\s_{K}|_{I}$. Then $\s_{I}$ is a splitting of $I$. Similarly, $\s_{I_{+}}:=\s_{K}|_{I_{+}}$ is a splitting of $I_{+}$. 

Since we are assuming that $\gcd(p,n)=1$, it follows that $I_{+}$ has a unique splitting. In contrast, the splitting of $I$ is not unique; every splitting of $I$ is of the form $\phi\cdot \s_{I}$, where $\phi\in\mathrm{Hom}(I,\mu_{n})\cong \mathrm{Hom}(T_{\kappa},\mu_{n})$. 

%Note that the character $\phi_{\s_{I}}|_{\s_{I}(I)}$ is by definition the trivial character. The pairs $(\It, \phi_{y}\cdot \phi_{\s_{I}})$, where $y\in Y$, are $\epsilon$-genuine depth $0$ minimal $\mathsf{K}$-types.

The next simple lemma highlights a special feature of the metaplectic Iwahori type theory. Recall that the important notion of associate types was introduced by Moy-Prasad \cite{MP94}. For the setting of tame covers see \cite[Definition 3.4]{HW09}. 

\begin{Prop}\label{Prop:Y_type_action}
    Let $y\in Y$ and $w\in W$. Let $\chi$ be any $\epsilon$-genuine depth $0$ character of $\It$.
    
    The two depth $0$ types $(\It,\chi)$ and $(\It,\phi_{y}\cdot \,^{w}\chi)$ are associate in $\Gt$.
\end{Prop}

\begin{proof}
    Let $\yt(\varpi)\in \mathrm{pr}^{-1}(y(\varpi))$ be arbitrary and let $k\in \Kt\cap N_{\Gt}(\Tt)$ be any representative of $w^{-1}\in W$. The result follows because for any $t\in \Tt_{0}$ the identity \eqref{Eq:Tt_Comm_Form} implies that $\chi(k\yt(\varpi)t\yt(\varpi)^{-1}k^{-1}) = (\phi_{y}\cdot \,^{w}\chi)(t)$.
\end{proof}

\begin{Rmk}
    From Lemma \ref{Lem:Char_equi} and Proposition \ref{Prop:Y_type_action} we get an action of the extended affine Weyl group $W\ltimes Y$, factoring through the finite group $W\ltimes Y/\Yt$, on the $\epsilon$-genuine depth $0$ Iwahori types. Moreover, the action of $Y/\Yt$ is free.
\end{Rmk}

The next corollary follows by combining Proposition \ref{Prop:Y_type_action} and specializing \cite[Proposition 3.7]{HW09}, which is the tame covering analog of \cite[Proposition 6.2]{MP96}.

\begin{Cor}\label{Cor:Translate_Types}
    Let $\pi\in\mathrm{Irr}_{gen}(\Gt)$. Let $\chi:\It\rightarrow \C^{\times}$ be a depth $0$ $\epsilon$-genuine character of $\It$.
    
    If $\pi^{(\It,\chi)}\neq 0$, then $\pi^{(\It,\phi_{y}\cdot \,^{w}\chi)}\neq 0$ for all $y\in Y$ and $w\in W$.
\end{Cor}

In the next proposition we describe how the appearance of the depth $0$ Iwahori type $(\s_{I}(I),\chi)$ and some root system combinatorics translates into the setting of finite groups of Lie type. Let $G_{\kappa}= K/K_{+}$, $B_{\kappa} = I/K_{+}$, and $U^{-}_{\kappa} = I^{-}_{+}/K_{+}$. Since $T_{\kappa}\cong I/I_{+}$, the natural quotient map gives a surjective group homomorphism $B_{\kappa}\rightarrow T_{\kappa}$.

\begin{Prop}\label{Prop:Finite_Generic}
    Let $\pi\in\mathrm{Irr}_{\epsilon}(\Gt)$. Let $\psi_{\kappa}$ be a non-degenerate character of $U^{-}_{\kappa}$. Let $\chi\in\mathrm{Hom}(T_{\kappa},\C^{\times})$, which can be inflated to a character of $\s_{I}(I)$. Let $\Ocal_{\chi}\subset \mathrm{Hom}(T_{\kappa},\C^{\times})$ be the $W\ltimes Y/\Yt$-orbit of $\chi$. 

    Let $N$ be equal to the number of free $W$-orbits in $\Ocal_{\chi}$.

    If $\pi^{(\s_{I}(I),\chi)}\neq 0$, then the $G_{\kappa}$-representation $\pi^{K_{+}}$ contains at least $N$ non-isomorphic $\psi_{\kappa}$-generic irreducible constituents. More specifically, these $N$ constituents can be taken to be regular principal series.
\end{Prop}

\begin{proof}
    By Proposition \ref{Cor:Translate_Types}, $\pi^{(\s_{I}(I),\chi^{\prime})}\neq 0$ for all $\chi^{\prime}\in \Ocal_{\chi}$. Let $\chi_{reg}\in \Ocal_{\chi}$ be any character such that the $W$-orbit of $\chi_{reg}$ is free, if it exists. Since the regular principal series $\mathrm{Ind}_{B_{\kappa}}^{G_{\kappa}}\chi_{reg}$ is irreducible, Frobenius reciprocity implies that  $\mathrm{Ind}_{B_{\kappa}}^{G_{\kappa}}\chi_{reg}$ is a $G_{\kappa}$-submodule of $\pi^{K_{+}}$. We also know that $\mathrm{Ind}_{B_{\kappa}}^{G_{\kappa}}\chi_{reg}$ is generic for any non-degenerate character. %Thus $\pi^{K_{+}}$ contains a $\psi_{\kappa}$-generic constituent. 

    Finally, since distinct free $W$-orbits in $\Ocal_{\chi}$ produce nonisomorphic regular principal series in $\pi^{K_{+}}$, it follows that the number of distinct $\psi_{\kappa}$-generic constituents is at least $N$.
\end{proof}

\section{Lifting Whittaker functionals}

In this section we prove Proposition \ref{Prop:Whit_lift}, which essentially relates Whittaker functionals on $\pi^{K_{+}}$ to Whittaker functionals on $\pi$. We suspect that this result is known, but we do not know a reference.

The sequence \eqref{Eqn:CExt} splits uniquely over any (algebraic) unipotent subgroup. We use this splitting to identify $U^{\pm}$ with its image in $\Gt$. Let $\psi\in\mathrm{Hom}(U^{-},\C^{\times})$ be a nondegenerate character of conductor $(\varpi)$. This means for every $\alpha\in\Delta$ the character $\psi$ is trivial on $U_{-\alpha,1}:=U_{-\alpha}\cap K_{+}$ and nontrivial on $U_{-\alpha,0}:=U_{-\alpha}\cap K$.

The inclusion $U^{-}\cap K\subset I^{-}_{+}$ descends to an isomorphism $(U^{-}\cap K)/(U^{-}\cap K_{+})\cong U^{-}_{\kappa}$. Since $\psi$ has conductor $(\varpi)$, the character $\psi|_{U^{-}\cap K}$ descends to a non-degenerate character of $(U^{-}\cap K)/(U^{-}\cap K_{+})$. Thus we may view $\psi|_{U^{-}\cap K}$ as a character of $U^{-}_{\kappa}$, which we call $\psi_{\kappa}$. We also write $\psi_{\kappa}$ for this character pulled back to $I^{-}_{+}$.

We begin with some notation. Let $\mathrm{Wh}_{\psi}(\pi):=\mathrm{Hom}_{G}(\pi,\mathrm{Ind}_{\mu_{n}U^{-}}^{G}\epsilon\otimes\psi)$ be the space of Whittaker models and $\mathcal{V}:=\mathrm{ind}_{\mu_{n}U^{-}}^{G}(\epsilon^{-1}\otimes\psi^{-1})$ the Gelfand-Graev representation. Let $\mathcal{W}:=\mathrm{ind}_{\mu_{n}I^{-}_{+}}^{G}\epsilon^{-1}\otimes\psi_{\kappa}^{-1}$. Let $\Hcal$ be the convolution algebra of locally constant compactly supported functions $f:\Gt\rightarrow \C$ such that $f(\zeta\gamma_{1}g\gamma_{2})=\epsilon(\zeta)^{-1}f(g)$, where $\zeta\in\mu_{n}$ and $\gamma_{j}\in \s_{I_{+}}(I_{+})$.

Recall that $\mathcal{W}^{\s_{I_{+}}(I_{+})}\cong \mathcal{V}^{\s_{I_{+}}(I_{+})}$ as $\Hcal$-modules, by \cite[Lemmas 5.1-5.3]{GGK24}. (This isomorphism was brought to our attention by the referee of \cite{GGK24}.) 

For $\pi$ a smooth admissible $\Gt$-representation, we write $\pi^{\vee}$ for its smooth dual.

\begin{Prop}\label{Prop:Whit_lift}
    Let $\pi\in\mathrm{Irr}_{\epsilon}(\Gt)$ such that $\pi^{\s_{I_{+}}(I_{+})}\neq 0$. Then as vector spaces
    \begin{equation*}
        \mathrm{Hom}_{\C}(\pi_{(U^{-},\psi)},\C)\cong [(\pi^{\vee})^{K_{+}}]^{(U^{-}_{\kappa},\psi_{\kappa}^{-1})}.
    \end{equation*}
\end{Prop}

\begin{proof}
    By Frobenius reciprocity and passing to coinvariants we have
    \begin{equation*}
        \mathrm{Hom}_{\C}(\pi_{(U^{-},\psi)},\C)\cong \mathrm{Wh}_{\psi}(\pi).
    \end{equation*}
    By dualizing we get
    \begin{equation*}
        \mathrm{Wh}_{\psi}(\pi)\cong \mathrm{Hom}_{G}(\mathcal{V},\pi^{\vee}).
    \end{equation*}
    Since $\pi^{\s_{I_{+}}(I_{+})}\neq 0$ we may use that the pair $(\s_{I_{+}}(I_{+}),\mathrm{triv})$ is a Bushnell-Kutzko type \cite[Lemma 5.8]{GGK24} to get
    \begin{equation*}
        \mathrm{Hom}_{G}(\mathcal{V},\pi^{\vee})\cong \mathrm{Hom}_{\Hcal}(\mathcal{V}^{\s_{I_{+}}(I_{+})},(\pi^{\vee})^{\s_{I_{+}}(I_{+})}).
    \end{equation*}
    By \cite[Lemmas 5.1-5.3]{GGK24}, it follows that $\mathcal{W}^{\s_{I_{+}}(I_{+})}\cong \mathcal{V}^{\s_{I_{+}}(I_{+})}$ as $\Hcal$-modules. Thus 
    \begin{equation*}
        \mathrm{Hom}_{\Hcal}(\mathcal{V}^{\s_{I_{+}}(I_{+})},(\pi^{\vee})^{\s_{I_{+}}(I_{+})})\cong \mathrm{Hom}_{\Hcal}(\mathcal{W}^{\s_{I_{+}}(I_{+})},(\pi^{\vee})^{\s_{I_{+}}(I_{+})}).
    \end{equation*}
    Again since $(\s_{I_{+}}(I_{+}),\mathrm{triv})$ is a Bushnell-Kutzko type we have 
    \begin{equation*}
        \mathrm{Hom}_{\Hcal}(\mathcal{W}^{\s_{I_{+}}(I_{+})},(\pi^{\vee})^{\s_{I_{+}}(I_{+})})\cong \mathrm{Hom}_{G}(\mathcal{W},\pi^{\vee}).
    \end{equation*}
    By Frobenius reciprocity and because $\psi$ is a one dimensional representation we have
    \begin{equation*}
        \mathrm{Hom}_{G}(\mathcal{W},\pi^{\vee})\cong \mathrm{Hom}_{I^{-}_{+}}(\psi_{\kappa}^{-1},\pi^{\vee})\cong (\pi^{\vee})^{(I^{-}_{+},\psi_{\kappa}^{-1})}.
    \end{equation*}
    Finally we can compute the space of invariants in two steps because
    \begin{equation*}
        (\pi^{\vee})^{(I^{-}_{+},\psi_{\kappa}^{-1})}\cong [(\pi^{\vee})^{K_{+}}]^{(U^{-}_{\kappa},\psi_{\kappa}^{-1})}.
    \end{equation*}
\end{proof}

\section{Main Theorem}

\begin{Thm}\label{Thm:generalized_main}
    Let $\pi\in \mathrm{Irr}_{\epsilon}(\Gt)$. Let $\chi\in\mathrm{Hom}(T_{\kappa},\C^{\times})$, which can be inflated to a character of $\s_{I}(I)$. Let $\Ocal_{\chi}\subset \mathrm{Hom}(T_{\kappa},\C^{\times})$ be the $W\ltimes Y/\Yt$-orbit of $\chi$. Let $N$ be the number of free $W$-orbits in $\mathcal{O}_{\chi}$. Let $\psi\in \mathrm{Hom}(U^{-},\C^{\times})$ be a non-degenerate character of conductor $(\varpi)$.

    If $\pi^{(\s_{I}(I),\chi)}\neq 0$, then the Whittaker dimension with respect to $\psi$ is at least $N$. 
\end{Thm}

\begin{proof}
    Since $\pi^{(\s_{I}(I),\chi)}\neq 0$, Proposition \ref{Prop:Finite_Generic} implies that $\pi^{K_{+}}$ as a $G_{\kappa}$-representation contains $N$ non-isomorphic regular principal series. Thus it follows that for any non-degenerate character $\psi:U^{-}\rightarrow \C^{\times}$ we have 
    \begin{equation*}
        \mathrm{dim}(\pi^{K_{+}})^{(U^{-}_{\kappa},\psi_{\kappa}^{-1})}\geq N.
    \end{equation*}
    Thus by Proposition \ref{Prop:Whit_lift},
    \begin{equation*}
        \mathrm{dim}\pi_{(U^{-},\psi)}\geq N.
    \end{equation*}
\end{proof}

Now Theorem \ref{Thm:Main_Thm} follows from Theorem \ref{Thm:generalized_main} by setting $\chi=\chi_{\mathrm{triv}}$ to be the trivial character. In this case, the map $\phi: Y/\Yt\rightarrow \Ocal_{\chi_{\mathrm{triv}}}$ defined by $y+\Yt\mapsto \phi_{y}$ is a bijection, because $\s_{I} = \s_{K}|_{I}$ is a $W$-equivariant splitting of $I$. By Lemma \ref{Lem:Char_equi}, $\phi$ is $W$-equivariant. Thus the number of free $W$-orbits in $Y/\Yt$ is equal to the number of free $W$-orbits in $\Ocal_{\chi_{\mathrm{triv}}}$.

%Next we introduce a genuine character of $\It$. We begin with a general statement. If $H\subset G$ is subgroup and $\s:H\rightarrow \overline{H}$ is a splitting of \eqref{Eqn:CExt} over $H$, then the map $\mu_{n}\times H\cong \overline{H}$ defined by $(\zeta,h)\mapsto \zeta \s(h)$ is a group isomorphism. 

%The map $\mu_{n}\times H\rightarrow \C^{\times}$ defined by $(\zeta,h)\mapsto \epsilon(\zeta)$ is a group homomorphism and if we pullback to $\overline{H}$ via the isomorphism $\mu_{n}\times H\cong \overline{H}$ we get an $\epsilon$-genuine character. We will write $\phi_{\s}:\overline{H}\rightarrow \C^{\times}$ for this $\epsilon$-genuine character.

\newpage

\begin{center}{\bf Acknowledgments}\end{center}

This work began with a conversation with Dipendra Prasad at the IASM-BIRS 5-day workshop on the wave front set and endoscopic classification in Hangzhou, China in October 2026. Prasad prompted the author with the question: Is there a group theoretic explanation for the genericity of genuine Iwahori-spherical representation of high degree BD covers?

I would also like to thank Wee Teck Gan and Gordan Savin for comments on an earlier draft.

\bibliographystyle{amsalpha}
\bibliography{Hecke_algebra2_bio}

\end{document}